\documentstyle[amssymb,10pt]{amsart}
\newtheorem{theorem}{Theorem}[section]
\newtheorem{lemma}[theorem]{Lemma}
\newtheorem{proposition}[theorem]{Proposition}
\newtheorem{corollary}[theorem]{Corollary} 
\theoremstyle{definition}

\newtheorem{remark}[theorem]{Remark}

\newcommand{\id}{\operatorname{id}} 
 
\newcommand{\Ker}{\text{Ker\,}}

\newcommand{\nc}{\newcommand}
\nc{\Symm}{{\on{Sym}}}

\newcommand{\on}{\operatorname}   
\newcommand{\eps}{\varepsilon}

 \nc{\cE}{{\cal E}}

\nc{\D}{{\mathfrak d}}
\nc{\SL}{{\mathfrak sl}}

\nc{\gt}{{\mathfrak gt}}
\nc{\grt}{{\mathfrak grt}}
\nc{\gtm}{{\mathfrak gtm}}
\nc{\grtm}{{\mathfrak grtm}}
\nc{\gtmd}{{\mathfrak gtmd}}
\nc{\grtmd}{{\mathfrak grtmd}}

\nc{\HH}{{\mathfrak h}}

\newcommand{\p}{{\mathfrak{p}}}

\renewcommand{\t}{{\mathfrak{t}}}

\newcommand{\f}{{\mathfrak{f}}}

\nc{\wh}{\widehat}\nc{\wt}{\widetilde}

\newcommand{\kk}{{\bf k}}

\newcommand{\Assoc}{{\bf Assoc}}

\newcommand{\ben}{\begin{enumerate}}
\newcommand{\een}{\end{enumerate}}

\newcommand{\CC}{{\mathbb{C}}}
\newcommand{\QQ}{{\mathbb{Q}}}

\renewcommand{\t}{{\mathfrak t}}

\begin{document}

\title[On the Drinfeld generators of $\grt_1(\kk)$ and $\Gamma$-functions for 
associators]{On the Drinfeld generators of $\grt_1(\kk)$ and 
$\Gamma$-functions for associators}

\begin{abstract} 
We prove that the Drinfeld generators of $\grt_1(\kk)$ span the 
image of this Lie algebra in the abelianization of the commutator 
of the free Lie algebra with two generators. We show that this 
result implies $\Gamma$-function formulas for arbitrary associators.  
\end{abstract}

\author{Benjamin Enriquez}
\address{IRMA (CNRS), rue Ren\'e Descartes, F-67084 Strasbourg, France}
\email{enriquez@@math.u-strasbg.fr}


\maketitle

\section*{Introduction and main results}

\subsection{Results on $\grt_1(\kk)$}

Let $A,B$ be free noncommutative variables and $\kk$ be a field
with $\on{char}(\kk) = 0$. Let $\f_2(A,B)$ be the free Lie algebra 
generated by $A,B$. The Lie algebra $\grt_1(\kk)$ is defined in \cite{Dr}
as the set of all $\psi\in \f_2(A,B)$, such that 
\begin{equation} \label{dual}
\psi(B,A) = -\psi(A,B), 
\end{equation}
\begin{equation} \label{hex}
\psi(A,B) + \psi(B,C) + \psi(C,A) = 0\ \on{if}
\ C = -A-B, 
\end{equation}
\begin{equation} \label{pent}
\psi^{12,3,4} - \psi^{1,23,4} + \psi^{1,2,34} = \psi^{2,3,4} + \psi^{1,2,3}.   
\end{equation}
The last relation takes place in the Lie algebra $\t_4$, 
defined as follows. When $n\geq 2$, $\t_n$ is the Lie algebra 
with generators $t_{ij}$, 
$i\neq j\in \{1,\ldots,n\}$ and relations $t_{ij} = t_{ji}$
if $i\neq j$, 
$[t_{ij} + t_{ik},t_{jk}] = 0$ if $i,j,k$ are distinct, and 
$[t_{ij},t_{kl}] = 0$ if $i,j,k,l$ are all distinct. 
If $I_1,\ldots,I_n$ are disjoint subsets of $\{1,\ldots,m\}$, then 
the Lie algebra morphism $\t_n \to \t_{m}$,    
$\psi \mapsto d^{I_1,\ldots,I_n}(\psi) = \psi^{I_1,\ldots,I_n}$ is defined by 
$t_{ij} \mapsto \sum_{\alpha\in I_i,\beta\in I_j} t_{\alpha\beta}$. 
Then $\t_3$ is the direct sum of its center $\kk(t_{12} + t_{13} + t_{23})$
and the free Lie algebra generated by $t_{12},t_{23}$, and we use the
identifications $A = t_{12}$, $B = t_{23}$. 

$\grt_1(\kk)$ is also equipped with a graded 
Lie algebra structure (it is not a Lie subalgebra of $\p$). 

Define $\p\subset \f_2(A,B)$ as the commutator subalgebra. 
If we assign degrees $1$ to $A$ and $B$, then $\p$ is the sum of 
all the components of $\f_2(A,B)$ of degree $>1$. 
Then $\p/[\p,\p]$ is an abelian Lie algebra, linearly spanned 
by the classes $p_{k\ell}$ of 
$\on{ad}(A)^{k-1} \on{ad}(B)^{\ell-1}([A,B])$, where $k,\ell\geq 1$. 
We have $\grt_1(\kk) \subset \p$; $\grt_1(\kk)$ is a graded 
subspace of $\p$. 

In \cite{Dr}, Drinfeld constructed a family of elements 
$\sigma_n \in \grt_1(\kk)$ ($n = 3,5,7,\ldots$), such that 
the image of the class $[\sigma_n]$ of $\sigma_n$
in $\p/[\p,\p]$
under the isomorphism $i : \p/[\p,\p] \simeq (\overline A\ \overline B) 
\subset \kk[\overline A,\overline B]$, 
$p_{k\ell}\mapsto \overline A^k \overline B^\ell$ is 
$$
i([\sigma_n]) = (\overline A+\overline B)^n -\overline A^n 
-\overline B^n. 
$$

We will prove: 

\begin{theorem} \label{main:thm}
Assume that $\psi\in\grt_1(\kk)$ is homogeneous of degree $n$. 
If $n$ is odd and $\geq 3$, then the 
image $[\psi]$ of $\psi$ in $\p/[\p,\p]$ is proportional 
to $[\sigma_n]$. Otherwise, this image is zero. 
\end{theorem}

\begin{remark} Actually, we will show that the conclusion of this 
theorem is valid if we assume that $\psi$ only satisfies 
(\ref{hex}) and (\ref{pent}). 
\end{remark}

\begin{remark} The maps $\psi\mapsto \psi^{12,3,4}$, 
$\psi\mapsto \psi^{1,23,4}$, etc., extend to algebra morphisms 
$U(\f_2(A,B)) \to U(\t_4)$. Similarly to \cite{EH}, one can show that 
$\{\psi\in \f_2(A,B) | \psi$ satisfies (\ref{dual}), (\ref{hex}) 
and (\ref{pent})$\} = \{\alpha^{12,3} - \alpha^{1,23} - \alpha^{2,3}
+ \alpha^{1,2} | \alpha\in t_{12} \kk[t_{12}]\} \oplus \grt_1(\kk)$. 
\end{remark}

\subsection{$\Gamma$-functions for associators}

Let $\wh F_2$ be the degree completion of $U(\f_2(A,B))$ ($A$ and $B$ 
have degree $1$). 

If $\lambda\in \kk^\times$, then $\Assoc_\lambda(\kk)$ is defined as 
the set of all $\Phi \in \wh F_2^\times$, such that 
$$
\Delta(\Phi) = \Phi \otimes \Phi, 
$$
$$
\Phi(A,B) e^{\lambda A/2} \Phi(C,A) e^{\lambda C/2} 
\Phi(B,C) e^{\lambda B/2} = 1\ \on{if}\ C = -A-B, 
$$
$$
\Phi^{2,3,4}\Phi^{1,23,4} \Phi^{1,2,3} = \Phi^{1,2,34} \Phi^{12,3,4}
$$ 
In particular, $\Phi$ has the form $\Phi = \on{exp}(\varphi)$, 
with $\varphi\in \wh\f_2(A,B)$ (the degree completion of 
$\f_2(A,B)$). 

We also set $\Assoc(\kk) = \{(\lambda,\Phi) | \lambda\in \kk^\times, 
\Phi\in \Assoc_\lambda(\kk)\}$. 

If $X$ is any element in $\wh F_2$, then there is a unique pair 
$(X_A,X_B)$ of elements of $\wh F_2$, such that 
$X = \eps(X)1 + X_A A + X_B B$ (here $\eps$ is the counit map of 
$\wh F_2$). We denote by $X\mapsto X^{\on{ab}}$ the abelianization 
morphism $\wh F_2 \to \kk[[\overline A,\overline B]]$, defined as 
the unique continuous algebra morphism such that $A\mapsto \overline A$, 
$B\mapsto \overline B$. 

Recall the formula $\zeta(n) = (2\pi i)^n r_n$ for $n$ even, 
where $r_n$ is a rational number (we have $r_n = -B_n /(2n!)$, 
where $B_n$ is the Bernoulli number defined by $u/(e^u-1) = 
\sum_{k\geq 0} B_k u^k/k!$). 

\begin{corollary} \label{cor:Gamma}
Let $\lambda\in\kk^\times$ and $\Phi\in \Assoc_\lambda(\kk)$, then 
there exists a unique sequence $(\zeta_\Phi(n))_{n\geq 2}$ of elements of 
$\kk$, such that 
\begin{equation} \label{Phi:Gamma}
(1 + \Phi_B B)^{\on{ab}} = {{\Gamma_\Phi(\overline A + \overline B)}
\over{\Gamma_\Phi(\overline A) \Gamma_\Phi(\overline B)}}, 
\end{equation}
where $\Gamma_\Phi$ is the invertible formal series 
$\Gamma_\Phi(u) = \on{exp}(-\sum_{n\geq 2} \zeta_\Phi(n) u^n/n)$. 
We have $\zeta_\Phi(n) = \lambda^n r_n$ for $n$ even. 
\end{corollary}

This result is contained in an unpublished paper by Deligne and 
Terasoma. Our proof relies on Theorem \ref{main:thm} and the torsor
structure of $\Assoc(\kk)$. 

\subsection*{Acknowledgements} I would like to thank G. Halbout for 
discussions on $\grt_1(\kk)$ in December 2002. I also thank G. Racinet
for informing me about the unpublished work of Deligne and Terasoma, 
and T. Terasoma for sending me a preliminary version of this work.

\section{Proof of Theorem \ref{main:thm}}

According to \cite{Dr}, the Lie algebras $\t_n$ have the following properties.   
The elements $t_{in}$, $i=1,\ldots,n-1$ generate a free 
subalgebra $\f_{n-1} \subset \t_n$. The Lie subalgebra of 
$\t_n$ generated by the $t_{ij}$, $i\neq j\in \{1,\ldots,n-1\}$
is isomorphic to $\t_{n-1}$. We have $\t_n = \f_{n-1} \oplus \t_{n-1}$; 
this is a semidirect product as $\t_{n-1}$ may be viewed as a Lie 
algebra of derivations of $\f_{n-1}$. 

Let us set $T_n = U(\t_n)$. The Lie algebra morphisms $\psi\mapsto 
d^{I_1,\ldots,I_n}(\psi) = \psi^{I_1,\ldots,I_n}$ extend to algebra 
morphisms $T_n \to T_m$, which we denote in the same way. 

We set $d = -d^{2,3,4} 
+ d^{12,3,4} - d^{1,23,4} + d^{1,2,34} - d^{1,2,3}$. 
So $d = d' + d''$, where $d' = -d^{2,3,4} + d^{12,3,4} - d^{1,23,4}$
and $d'' = d^{1,2,34} - d^{1,2,3}$. 
Then $d,d',d''$ are linear maps $T_3 \to T_4$, which restrict to
linear maps $\t_3\to \t_4$ (which we denote the same way). 

\begin{lemma}
The linear maps $d,d'$ and $d''$ map $\f_2 \subset \t_3$ to 
$\f_3 \subset \t_4$. 
\end{lemma} 

{\em Proof.} There is a unique Lie algebra morphism 
$\eps_4 : \t_4 \to \t_3$, with $\eps_4(t_{ij}) = t_{ij}$ for $i<j<4$
and $\eps_4(t_{i4}) = 0$ for $i<4$. Then $\on{Ker}(\eps_4) = \f_3$. 

We have $\eps_4 \circ d'' = 0$, hence $d''(\f_2) \subset \on{Ker}(\eps_4) 
= \f_3$. On the other hand, 
the Lie algebra morphisms $d^{2,3,4},d^{12,3,4}$ and 
$d^{1,23,4} : \t_3 \to \t_4$ are such that 
$(t_{13},t_{23}) \mapsto (t_{24},t_{34})$, 
$(t_{13},t_{23}) \mapsto (t_{14} + t_{24},t_{34})$, 
$(t_{13},t_{23}) \mapsto (t_{14}, t_{24}+t_{34})$, 
so these morphisms take the generators of $\f_2$ to $\f_3$, 
so they induce Lie algebra morphisms $\f_2 \to \f_3$. 
Therefore $d'(\f_2) \subset \f_3$. It follows that 
$d(\f_2) \subset \f_3$. 
\hfill \qed \medskip 

We set $F_{n-1} := U(\f_{n-1})$, $T_n := U(\t_n)$. Then the tensor 
product of inclusions followed by multiplication induces a linear 
isomorphism $F_{n-1} \otimes T_{n-1} \stackrel{\simeq}{\to} T_n$.  
We denote by $\Pi : T_4 \to F_4$ the composition 
$F_4 \stackrel{\simeq}{\to} F_3 \otimes T_3 \stackrel{\id\otimes \eps}{\to}
F_3$, where $\eps : T_3 = U(\t_3) \to \kk$ is the counit map. 

\begin{lemma}
$d' : T_3 \to T_4$ is such that $d'(F_2) \subset F_3$. 
On the other hand, the composition $F_2 \stackrel{d''}{\to} T_3 
\stackrel{\Pi}{\to} F_3$ is a linear map 
$F_2\to F_3$ extending $d'' : \f_2\to \f_3$. So we have commuting 
diagrams 
$$
\matrix 
F_2 & \stackrel{d'}{\to} & F_3 \\ 
\cup & & \cup \\ 
\f_2 & \stackrel{d'}{\to} & \f_3
\endmatrix
\quad {\rm and} \quad 
\matrix 
F_2 & \stackrel{d''}{\to} & T_3 & \stackrel{\Pi}{\to}& F_3 \\ 
\cup & & & & \cup \\ 
\f_2 & & \stackrel{d''}{\to} & & \f_3
\endmatrix
$$
\end{lemma}

{\em Proof.} We have seen that the Lie algebra morphisms 
$d^{2,3,4},d^{12,3,4}$ and $d^{1,23,4} : \t_3\to \t_4$
restrict to Lie algebra morphisms $\f_2\to \f_3$. 
It follows that their extensions to algebra morphisms 
$T_3\to T_4$ restrict to algebra morphisms $F_2 \to F_3$. 
As $d'$ is a linear combination of these morphisms, it follows 
that $d'(F_2) \subset F_3$. 

Let $\psi$ be an element of $\f_2$. We have seen that 
$d''(\psi)\in \f_3 \subset F_3$. For any $x\in F_3$, we have 
$\Pi(x)=x$. Therefore $\Pi (d''(\psi)) = d''(\psi)$. 
\hfill \qed \medskip 

It follows that 
\begin{equation} \label{lanzmann}
\on{Ker}(d : \f_2 \to \f_3) = \on{Ker} (d' + \Pi \circ d'' : 
F_2 \to F_3) \cap \f_2.
\end{equation} 

We now define vector subspaces $I_2 \subset F_2$ and $I_3 \subset F_3$
as follows. 

Set $X := t_{13}$, $Y:= t_{23}$ (elements of $F_2$). 
Then $F_2 = \kk\langle X,Y \rangle$. A basis of $F_2$ is the 
set of all words in $X,Y$. We define $I_2$ to be the linear 
span of all words of the form $wXw'Yw''$, where $w,w',w''$ are 
words in $X,Y$. So $I_2$ is 
spanned by the non-lexicographically ordered words, where the order 
is $Y<X$. 

Set $x := t_{14}$, $y :=t_{24}$, $z :=t_{34}$ (elements of $F_3$). 
Then $F_3 = \kk\langle x,y,z \rangle$. A basis of $F_3$ is the 
set of all words in $x,y,z$. We define $I_3$ to be the linear 
span of all words of the form $wxw'yw''$, $wxw'zw''$ or 
$wyw'zw''$, where $w,w',w''$ are words in $x,y,z$. So $I_3$ is 
spanned by the non-lexicographically ordered words, where the order 
is $z<y<x$. 

\begin{lemma} \label{lemma:quot}
We have linear isomorphisms $F_2/I_2 \simeq \kk[\overline X,\overline Y]$
and $F_3/I_3 \simeq \kk[\overline x,\overline y,\overline z]$, 
where $\overline X,\overline Y$ on one hand, $\overline x,\overline y,
\overline z$ on the other hand are free commutative variables. 
\end{lemma}

{\em Proof.} Let $V_2 \subset F_2$ be the subspace with basis 
$Y^bX^a$, where $a,b\geq 0$. Then $V_2 \oplus I_2 = F_2$, 
so we have an isomorphism $F_2/I_2 \simeq V_2$. We then 
compose this isomorphism with $V_2 \to \kk[\overline X,\overline Y]$, 
$Y^bX^a \mapsto \overline X^a \overline Y^b$. 

Let $V_3 \subset F_3$ be the subspace with basis 
$z^cy^bx^a$, where $a,b,c\geq 0$. Then $V_3 \oplus I_3 = F_3$, 
so we have an isomorphism $F_3/I_3 \simeq V_3$. We then 
compose this isomorphism with $V_3 \to \kk[\overline x,\overline y,
\overline z]$, $z^cy^bx^a \mapsto \overline x^a \overline y^b\overline z^c$. 
\hfill \qed \medskip 

\begin{remark}
Even though $I_\alpha$ is a two-sided ideal of $F_\alpha$
($\alpha = 2,3$), the isomorphisms of Lemma \ref{lemma:quot}
are not algebra isomorphisms. Indeed, the algebras $F_\alpha/I_\alpha$
are noncommutative and have zero divisors. \hfill \qed \medskip 
\end{remark}

\begin{lemma}
Define $\pi : F_2\to F_3$ are the composition 
$F_2 \stackrel{d^{1,2,34}}{\to} T_3 \stackrel{\Pi}{\to} F_3$. 
Then $\Pi\circ d''(a) = \pi(a) - \eps(a)1$ for any 
$a\in F_2$ (here $\eps : F_2 = U(\f_2) \to \kk$ is the counit map). 

Let $\tau_{13},\tau_{23}$ be the derivations of $F_3$ defined by 
$\tau_{13}: x\mapsto [x,z]$, $y\mapsto 0$, $z\mapsto [z,x]$
and $\tau_{23} : x\mapsto 0$, $y\mapsto [y,z]$, $z\mapsto [z,y]$.  

Then we have, for any $a\in F_2$, 
\begin{equation} \label{induction}
\pi(Xa) = x\pi(a) + \tau_{13}(\pi(a)), \quad 
\pi(Ya) = y\pi(a) + \tau_{23}(\pi(a)).  
\end{equation}
\end{lemma}

{\em Proof.}
If $a\in T_3$, then $\Pi\circ d^{1,2,3}(a) = \eps(a)1$, 
where $\eps : T_3 = U(\t_3) \to \kk$ is the counit map. 
So if $a\in F_2$, we have $\Pi\circ d''(a) = 
\Pi\circ d^{1,2,34}(a) - \eps(a)1$. 

Let us now prove formulas (\ref{induction}).  
Let $a\in F_2$, then $d^{1,2,34}(a) = \pi(a) 
+ \sum_i a_it_i$, where $a_i\in F_3$ and $t_i
\in \Ker(\eps : T_3 = U(\t_3) \to \kk)$ (here $\eps$
is the counit map). So 
\begin{align*}
& d^{1,2,34}(a) = (t_{13} + t_{14})
(\pi(a) + \sum_i a_i t_i)
\\ &  = x (\pi(a) + \sum_i a_i t_i) + [t_{13},\pi(a)] 
 + \pi(a) t_{13}
 + \sum_i [t_{13},a_i]t_i + \sum_i a_i (t_{13} t_i). 
\end{align*}
We have $[t_{13},b] = \tau_{13}(b)$ for any $b\in F_2$, 
so this is the sum of $x\pi (a) + \tau_{13}(\pi(a))$
and $\sum_i xa_i t_i + \pi(a) t_{13} + \sum_i \tau_{13}(a_i) t_i
+ \sum_i a_i (t_{13}t_i)$. The first term belongs to 
$F_3$ and the second term belongs to $F_3 \Ker(\eps : T_3\to\kk)$, 
so the image of their sum by $\Pi$ is the first term,
i.e.,  $x\pi (a) + \tau_{13}(\pi(a))$. This proves the first identity 
of (\ref{induction}). The second identity is proved in the same way. 
\hfill \qed \medskip 

\begin{proposition}
We have $d'(I_2) \subset I_3$ and $\Pi\circ d''(I_2) \subset I_3$. 
\end{proposition}

{\em Proof.} Let $w,w',w''$ are words in $X,Y$. Then 
$d^{2,3,4}(wXw'Yw'') = d^{2,3,4}(w) y d^{2,3,4}(w') z d^{2,3,4}(w'')$, 
which decomposes as a sum of words of the form $w_1 y w_2 z w_3$, 
where $w_i$ are words on $x,y,z$. So $d^{2,3,4}(wXw'Yw'') \in I_3$, 
which shows that $d^{2,3,4}(I_2) \subset I_3$. 

In the same way, $d^{12,3,4}(wXw'Yw'') = d^{12,3,4}(w) (x+y) 
d^{12,3,4}(w') z d^{12,3,4}(w'')$ belongs to $I_3$, so  
$d^{12,3,4}(I_2) \subset I_3$. 

We also have $d^{1,23,4}(wXw'Yw'') = d^{12,3,4}(w) x 
d^{12,3,4}(w') (y+z) d^{12,3,4}(w'')$ belongs to $I_3$, so  
$d^{1,23,4}(I_2) \subset I_3$. 

Since $d' = - d^{2,3,4} + d^{12,3,4} - d^{1,23,4}$, it follows 
that $d'(I_2) \subset I_3$. 

Let us now prove that $\Pi\circ d''(I_2) \subset I_3$.

If $a\in F_2$, we have $\Pi\circ d''(a) = 
\pi(a) - \eps(a)1$. Now $\eps(a)=0$
if $a\in I_2$, so we have to prove that 
$\pi(I_2) \subset I_3$.  

We denote by $I_2[n]$ the degree $n$ part of $I_2$
(here $X,Y$ have degree $1$). We will prove by induction on 
$n$ that $\pi(I_2[n]) \subset I_3$.  

When $n=2$, $I_2[n]$ is spanned by $XY$. Then 
$\pi(XY) = \Pi((t_{13} + x)(t_{23} + y)) = xy + \tau_{13}(y) 
= xy$ belongs to $I_3$. 

Let $n\geq 3$ and assume that we have proved that 
$\pi(I_2[n-1]) \subset I_3$.  Let us prove that 
$\pi(I_2[n]) \subset I_3$.  

$I_2[n]$ is spanned by the words $Yw X w' Y w''$, where 
$w,w',w''$ are words in $X,Y$ of total length $n-3$, and 
by the words $XY^bX^a$, where $a\geq 0,b>0$ and $a+b=n-1$.
We should prove that $\pi$ takes these words to $I_3$.   

Let us study the image of the first family of words. 
Set $w''' := wXw'Yw''$. Then $w'''\in I_2[n-1]$, and
according to (\ref{induction}),  
\begin{equation} \label{weiss}
\pi(Yw''') = y \pi(w''') + \tau_{23}(\pi(w''')).
\end{equation} 
Now the induction hypothesis implies that $\pi(w''') \in I_3$. 
Since $I_3$ is a two-sided ideal of $I_3$, $y\pi(w''')\in I_3$. 

On the other hand, let us prove that if $a\in I_3$, then 
$\tau_{23}(a)\in I_3$. If $a$ is a word of the form 
$wxw'yw''$, then 
$$   
\tau_{23}(a) = \tau_{23}(w)xw'yw'' + wx\tau_{23}(w')yw''
+ wxw'[y,z]w'' \in I_3. 
$$
If $a$ has the form $wxw'zw''$, then 
$$   
\tau_{23}(a) = \tau_{23}(w)xw'zw'' + wx\tau_{23}(w')zw''
+ wxw'[z,y]w'' \in I_3. 
$$
If $a$ has the form $wyw'zw''$, then 
$$
\tau_{23}(a) = \tau_{23}(w)yw'zw'' + w[y,z]w'zw'' 
+ wy\tau_{23}(w')zw'' + wyw'[z,y]w'' + wyw'z\tau_{23}(w'') \in I_3. 
$$
By linearity, it follows that if $a\in I_3$, then 
$\tau_{23}(a)\in I_3$. 

In particular, $\tau_{23}(\pi(w''')) \in I_3$. 
Therefore (\ref{weiss}) implies that $\pi(Yw''')\in I_3$. 

Let us now study the image of the second family of words, 
i.e. $\pi(XY^bX^a)$, where $b>0$ and $a+b=n-1$.

Let us first show: 

\begin{lemma} \label{lemma:1}
$\pi(Y^bX^a)$ has positive 
valuation\footnote{The valuation in $x_i$ of a nonzero 
element of a free algebra $k\langle x_1,\ldots,x_n\rangle$ 
is the smallest degree in $x_i$ of a word appearing with 
a nontrivial coefficient in its decomposition; the valuation 
of $0$ is $+\infty$.} 
in $y$. 
\end{lemma}

{\em Proof of Lemma.} Let 
$\eps_y : F_3 \to F_2 = \kk\langle t_{13},t_{23}\rangle$ 
be the morphism defined 
by $x\mapsto t_{13}$, $y\mapsto 0$, $z\mapsto t_{23}$. We want to 
show that $\eps_y \circ \pi(Y^bX^a)=0$.  

Recall that $\t_3 = \f_2 \oplus \t_2$, therefore the 
composed map $F_2 \otimes T_2 \to T_3^{\otimes 2} \to T_3$
(the first map is the tensor product of inclusions, the second 
map is the product) is a linear isomorphism $F_2\otimes T_2 
\stackrel{\simeq}{\to} T_3$. We denote by $\Pi_3$ the composed map 
$T_3 \stackrel{\simeq}{\to} F_2 \otimes T_2 
\stackrel{\id\otimes \eps}{\to} F_2$. 
Then we have a commutative diagram 
$$
\matrix 
T_4 & \stackrel{\eps_2}{\to} & T_3 \\ 
\scriptstyle{\Pi}\downarrow & & \downarrow\scriptstyle{\Pi_3}\\ 
F_3 & \stackrel{\eps_y}{\to} & F_2 
\endmatrix
$$
Here $\eps_2 : T_4\to T_3$ is the morphism induced by the 
Lie algebra morphism $\t_4\to\t_3$, $t_{ij}\mapsto 0$ if $i$ or $j=2$, 
$t_{13}\mapsto t_{12}$, $t_{14}\mapsto t_{13}$, $t_{34} \mapsto t_{23}$. 
Indeed, 
$$
P(t_{14},t_{24},t_{34})Q(t_{12},t_{13},t_{23})
\stackrel{\eps_2}{\to} P(t_{13},0,t_{23}) Q(0,t_{12},0) 
\stackrel{\Pi_3}{\to} P(t_{13},0,t_{23}) Q(0,0,0),
$$ whereas
$$
P(t_{14},t_{24},t_{34})Q(t_{12},t_{13},t_{23})
\stackrel{\Pi}{\to} P(t_{14},t_{24},t_{34}) Q(0,0,0) 
\stackrel{\eps_y}{\to} P(t_{13},0,t_{23}) Q(0,0,0).
$$ 

It follows that $\eps_y \circ \pi(Y^bX^a) = \Pi_2\circ \eps_2 \circ 
d^{1,2,34}(Y^bX^a)$. Now $\eps_2\circ d^{1,2,34}(Y^bX^a)
= (t_{23} + t_{24})^b (t_{13} + t_{14})^a = 0$ since $b>0$. 
\hfill \qed \medskip 

\begin{lemma} \label{lemma:2}
If $w$ is a word in $x,y,z$ of
positive degree in $y$, then 
\begin{equation} \label{partial}
xw + \tau_{13}(w) \in I_3.
\end{equation} 
\end{lemma}

{\em Proof of Lemma.}
The word $xw$ contains $xy$ as a subword, hence $xw\in I_3$. 
Let us now write $w$ as a product $w'yw''$, where $w',w''$
are words. We have $\tau_{13}(w) = \tau_{13}(w')yw'' + 
w'y\tau_{13}(w'')$. In general, if $w'''$ is a word, 
then $\tau_{13}(w''')$ has positive valuation both in 
$x$ and $z$. Since $\tau_{13}(w')$ (resp., $\tau_{13}(w'')$) 
has positive valuation in $x$ (resp., in $z$), 
$\tau_{13}(w')yw''$ (resp., $w'y\tau_{13}(w'')$) contains $xy$ 
(resp., $yz$) as a subword. If follows that $\tau_{13}(w) \in I_3$. 
This implies (\ref{partial}).  \hfill \qed \medskip 

{\em End of proof of Proposition.}
Now (\ref{induction}), Lemma \ref{lemma:1} and Lemma \ref{lemma:2}
imply that $\pi(XY^bX^a)\in I_3$. 
\hfill \qed \medskip 

It follows that $d'$ and $\Pi \circ d''$ induce maps 
$F_2/I_2 \to F_3/I_3$, which we compute explicitly. 

\begin{lemma}
The maps $\overline d', \overline d'' : 
\kk[\overline X,\overline Y] \to \kk[\overline x,\overline y,
\overline z]$ induced by $d'$ and $\Pi \circ d''$ are 
given by $\overline d' = -\overline d^{2,3,4} 
+ \overline d^{12,3,4} - \overline d^{1,23,4}$, where 
$$
\overline d^{2,3,4} : f(\overline X,\overline Y) \mapsto 
f(\overline y,\overline z), 
$$
$$
\overline d^{12,3,4} : f(\overline X,\overline Y) \mapsto 
{{\overline x f(\overline x,\overline z)
- \overline y f(\overline y,\overline z)}\over 
{\overline x - \overline y}},  
$$
$$
\overline d^{1,23,4} : f(\overline X,\overline Y) \mapsto 
{{\overline y f(\overline x,\overline y)
- \overline z f(\overline x,\overline z)}\over 
{\overline y - \overline z}},    
$$
and by 
$$
\overline d'' : 1 \mapsto 0, \; 
\overline X f(\overline X) \mapsto \overline x f(\overline x - \overline z), 
\;
\overline Y f(\overline Y) \mapsto \overline y f(\overline y - \overline z), 
\;
\overline X\ \overline Yf(\overline X,\overline Y)\mapsto 
\overline x\ \overline yf(\overline x-\overline z,\overline y-\overline z). 
$$
Here $f(\overline X),f(\overline Y)$ (resp., $f(\overline X,\overline Y)$)
are arbitrary $1$-variable (resp., $2$-variable) polynomials. 
\end{lemma}

{\em Proof.}
The maps $d^{2,3,4},d^{12,3,4}$ and $d^{1,23,4}$ all take 
$I_2$ to $I_3$, so they induce maps $\overline d^{2,3,4}, 
\overline d^{12,3,4}$ and $\overline d^{1,23,4}:F_2/I_2 \to F_3/I_3$. 
We then have $\overline d' = - \overline d^{2,3,4} + \overline d^{12,3,4} -
\overline d^{1,23,4}$. 

Let us compute $\overline d^{2,3,4}$. We have $d^{2,3,4}(Y^b X^a) = 
z^b y^a$, whose image in $F_3/I_3$ is $\overline z^b \overline y^a$. 
This implies the formula for $\overline d^{2,3,4}$. 

Let us compute $\overline d^{12,3,4}$. We have $d^{12,3,4}(Y^bX^a)
= z^b(x+y)^a$, whose projection on $V_3$ along $I_3$  
is $z^b(y^a + y^{a-1}x + \cdots + x^a)$. The image of 
this element in $\kk[\overline x,\overline y,\overline z]$ is
$$
(\overline x^a + \overline x^{a-1}\overline y + \cdots + \overline y^a)
\overline z^b = {{\overline x^{a+1} - \overline y^{a+1}}\over
{\overline x - \overline y}} \overline z^b.
$$ 
The formula for $\overline d^{12,3,4}$ follows by linearity.  

One computes $\overline d^{1,23,4}$ in the same way. 
We have $d^{1,23,4}(Y^bX^a)
= (y+z)^b x^a$, whose projection on $V_3$ along $I_3$ is 
is $(z^b + z^{b-1}y + \cdots + y^b) x^a$. The image of 
this element in $\kk[\overline x,\overline y,\overline z]$ is
$$
\overline x^a
(\overline y^b + \overline y^{b-1}\overline z + \cdots + \overline z^b)
= \overline x^a {{\overline y^{b+1} - \overline z^{b+1}}\over
{\overline y - \overline z}}.
$$ 
The formula for $\overline d^{1,23,4}$ follows by linearity.  

Let us now compute $\overline d''$. Clearly $\overline d''(1)=0$. 
Let $n>0$ and let us compute $\overline d''(\overline X^n)$. 
This is the image in $F_3/I_3$ of $\Pi \circ d^{1,2,34}(X^n)
= \Pi((t_{13} + t_{14})^n)$. We have  
$$
(t_{13} + t_{14})^n = ((t_{14} + t_{34} + t_{13}) - t_{34})^n 
= \sum_{k=0}^n (-1)^k 
C_n^k (t_{34})^k (t_{14} + t_{34} + t_{13})^{n-k} , 
$$
since $[t_{14} + t_{34} + t_{13},t_{34}]=0$. 
Now $[t_{13},t_{14}+t_{34}] = 0$, hence 
$$
(t_{34})^k (t_{14} + t_{34} + t_{13})^{n-k} 
= \sum_{\alpha=0}^{n-k} C_{n-k}^\alpha 
(t_{34})^k (t_{14} + t_{34})^{n-k-\alpha} (t_{13})^\alpha, 
$$ 
which is mapped by $\Pi$ to $(t_{34})^k (t_{14} + t_{34})^{n-k} 
= z^k (x+z)^{n-k}$. The projection of this element on $V_3$ along $I_3$
is $z^k(z^{n-k} + z^{n-k-1}x + \cdots + x^{n-k})$, whose image in 
$\kk[\overline x,\overline y,\overline z]$ is 
$\overline z^k (\overline x^{n-k+1} - \overline z^{n-k+1})/(\overline x -
\overline z)$. So 
$$
\overline d''(\overline X^n) = \sum_{k=0}^n (-1)^k C_n^k 
{{\overline x^{n-k+1} - \overline z^{n-k+1}}\over
{\overline x - \overline z}}
\overline z^k = \overline x(\overline x - \overline z)^{n-1}. 
$$
This implies the formula for $\overline d''(Xf(X))$. The formula for 
$\overline d''(\overline Yf(\overline Y))$ is proved in the same way. 

Let us now prove by induction on $k+\ell$ that when $k,\ell>0$,  
\begin{equation} \label{induct}
\overline d''(\overline X^k \overline Y^\ell) = \overline x\ \overline y
(\overline x - \overline z)^{k-1} (\overline y - \overline z)^{\ell-1}. 
\end{equation}

When $k=\ell=1$, $d^{1,2,34}(YX) = (t_{23}t_{13})^{1,2,34} 
\stackrel{\Pi}{\to} t_{24} t_{14} + [t_{23},t_{14}] = yx$ hence 
$\overline d''(\overline X\ \overline Y) = \overline x\ \overline y$, 
which proves (\ref{induct}) in this case.  

Assume that (\ref{induct}) holds for $k+\ell<n$ and let us prove it for 
$k+\ell=n$ ($k,\ell>0$). 

When $\ell=1$, we have $\overline d''(\overline X^k) = 
\overline x (\overline x - \overline z)^{k-1}$, therefore 
\begin{equation} \label{form:pi:Xk}
\pi(X^k) = \sum_{\alpha = 0}^{k-1} (-1)^\alpha C_{k-1}^\alpha z^\alpha
x^{k-\alpha} + \xi, \quad \on{where}\quad \xi\in I_3. 
\end{equation}
Then (\ref{induction}) implies that 
$$
\pi(YX^k) = y\pi(X^k) + \sum_{\alpha=0}^{k-1} (-1)^\alpha C_{k-1}^\alpha
\tau_{23}(z^\alpha x^{k-\alpha}) + \tau_{23}(\xi). 
$$
Now (\ref{form:pi:Xk}) implies that $y\pi(X^k)\in yx^k + I_3$. 
The projection of
$\tau_{23}(z^\alpha x^{k-\alpha})$ on $V_3$ along $I_3$ is $z^\alpha y
x^{k-\alpha}$ if $\alpha\neq 0$, $0$ otherwise. 

\begin{lemma} \label{lemma:tau:23}
$\tau_{23}(I_3) \subset I_3$. 
\end{lemma} 

{\em Proof of Lemma.}  If $w,w'$ and $w''$ are any words in $x,y,z$, then 
$$
\tau_{23}(wxw'yw'') = \tau_{23}(w)xw'yw'' + wx\tau_{23}(w')yw'' 
+ wxw'[y,z]w'' + wxw'y\tau_{23}(w'') \in I_3. 
$$
One proves similarly that $\tau_{23}(wxw'zw'')$ and $\tau_{23}(wyw'zw'')$
belong to $I_3$. 
\hfill \qed \medskip 

So $\overline d''(\overline X^k \overline Y) = 
\overline x^k \overline y + \sum_{\alpha=1}^{k-1} (-1)^\alpha
C_{k-1}^\alpha \overline x^{k-\alpha} \overline y\ \overline z^\alpha = 
\overline x\ \overline y (\overline x - \overline z)^{k-1}$, 
which proves (\ref{induct}) in this case. 

When $\ell>1$, we use (\ref{induct}) for $(k,\ell-1)$.
This gives 
$$
\pi(Y^{\ell-1} X^k) = \sum_{\alpha=0}^{k=1} \sum_{\beta=0}^{\ell-2}
(-1)^{\alpha+\beta} C_{k-1}^\alpha C_{\ell-1}^\beta z^{\alpha+\beta}
y^{\ell-1-\beta} x^{k-\alpha} + \eta,  
$$
where $\eta\in I_3$. 

Then (\ref{induction}) implies that 
$$
\pi(Y^\ell X^k) = y\pi(Y^{\ell-1}X^k) + \sum_{\alpha=0}^{k-1}
\sum_{\beta=0}^{\ell-2} (-1)^{\alpha+\beta} C_{k-1}^\alpha C_{\ell-1}^\beta
\tau_{23}(z^{\alpha+\beta} y^{\ell-1-\beta} x^{k-\alpha}) + \tau_{23}(\eta). 
$$
All the terms  in the expansion of $y\pi(Y^{\ell-1} X^k)$ belong to $I_3$,
except the terms corresponding to $\alpha = \beta =0$, so 
$y\pi(Y^{\ell-1} X^k)\in y^\ell x^k + I_3$.  

If $a,b,c\geq 0$, then the projection of $\tau_{23}(z^c y^b x^a)$ on 
$V_3$ along $I_3$ is $z^c (y-z) y^b x^a$ if $b\neq 0$ and $c\neq 0$; 
it is $z^c y^{b+1} x^a$ if $c\neq 0$ and $b=0$; it is 
$-z^{c+1}y^bx^a$ if $c=0$ and $b\neq 0$; and it is $0$ if $b=c=0$. 

Lemma \ref{lemma:tau:23} implies that $\tau_{23}(\eta)\in I_3$. Then the
projection of $\pi(Y^\ell X^k)$ on $V_3$ along $I_3$ is 
\begin{align*}
& y^\ell x^k - z y^{\ell-1} x^k + \sum_{(\alpha,\beta) \in (\{0,\ldots,k-1\} 
\times \{0,\ldots,\ell-2\}) - \{(0,0)\}} 
(-1)^{\alpha+\beta} C_{k-1}^\alpha C_{\ell-2}^\beta
z^{\alpha+\beta} (y-z) y^{\ell-1-\beta} x^{k-\alpha} 
\\ & = 
\sum_{\alpha=0}^{k-1} \sum_{\beta=0}^{\ell-1}
(-1)^{\alpha+\beta} C_{k-1}^\alpha C_{\ell-2}^\beta
z^{\alpha+\beta} (y-z) y^{\ell-1-\beta} x^{k-\alpha} . 
\end{align*}
So $\overline d''(\overline X^k \overline Y^\ell) = 
(\overline y - \overline z) \overline x\ \overline y
(\overline x - \overline z)^{k-1} (\overline y - \overline z)^{\ell-2}
= \overline x\ \overline y
(\overline x - \overline z)^{k-1} (\overline y - \overline z)^{\ell-1}$, 
which proves (\ref{induct}) in this case. This proves the induction. 
\hfill \qed \medskip 

\begin{lemma} \label{ker:bar:d}
Define $\overline d = \overline d' + \overline d'' : 
\kk[\overline X,\overline Y] \to 
\kk[\overline x,\overline y,\overline z]$. The kernel of 
$\overline d_{|(\overline X\ \overline Y)} : (\overline X\ \overline Y)
\to \kk[\overline x, \overline y,\overline z]$ is equal to the linear span of
the $\overline X \big( \overline X^n - \overline Y^n 
- (\overline X - \overline Y)^n \big) / (\overline X - \overline Y)$, 
where $n\geq 2$. 
\end{lemma}

{\em Proof.}
The map $\overline d_{|(\overline X\ \overline Y)}$ takes 
$f(\overline X,\overline Y)$ to 
\begin{align*}
& -f(\overline y,\overline z) 
+ {{\overline x f(\overline x, \overline z)
- \overline y f(\overline y,\overline z)}\over{\overline x - \overline y}}
- {{\overline y f(\overline x, \overline y)
- \overline z f(\overline x,\overline z)}\over{\overline y - \overline z}}
+ {{\overline x\ \overline y 
f(\overline x-\overline z,\overline y-\overline z)}\over
{(\overline x-\overline z)(\overline y-\overline z)}}
\\ & 
= 
{{\overline x\ \overline y}\over
{(\overline x- \overline y)(\overline y-\overline z)}}
\big( g(\overline x - \overline z,\overline y - \overline z) 
+ g(\overline x,\overline z) - g(\overline y,\overline z) 
- g(\overline x,\overline y) \big),  
\end{align*}
where $g(\overline x,\overline y) = 
{{\overline x - \overline y}\over {\overline x}} f(\overline x,\overline y)$. 

So $f(\overline X,\overline Y) \in \Ker(\overline d) \cap (\overline X\
\overline Y)$ iff $f(\overline X,\overline Y) \in (\overline X\ \overline Y)$
and 
\begin{equation} \label{eq:g}
g(\overline x - \overline z,\overline y - \overline z) 
+ g(\overline x,\overline z) - g(\overline y,\overline z) 
- g(\overline x,\overline y) = 0. 
\end{equation}

Let us solve (\ref{eq:g}), where $g(\overline x,\overline y) \in 
\kk[\overline x,\overline y]$. By linearity, we may assume that $g$ 
is homogeneous; let $n$ be its degree. If $n=0$, we get $g =$ a 
constant polynomial. Assume that $n>0$. Applying 
$(\partial / \partial\overline z)_{|\overline z = 0}$ to (\ref{eq:g}), 
we get 
$$
({\partial \over{\partial\overline x}} + {\partial \over{\partial\overline y}})
g(\overline x,\overline y) = c(\overline x^{n-1} - \overline y^{n-1}), 
$$
for some $c\in\kk$. This gives $c(\overline X,\overline Y) = 
h(\overline X - \overline Y) + c(\overline X^n - \overline Y^n)/n$, 
there $h(\overline X)\in \kk[\overline X]$ has degree $n$, so 
$g(\overline X,\overline Y) = c(\overline X^n - \overline Y^n)/n + 
c'(\overline X - \overline Y)^n$, for some $c'\in\kk$. 

Substituting this is (\ref{eq:g}), we get $c' = -c/n$, so the set of solutions 
of degree $n$ of (\ref{eq:g}) is the linear span of 
$g(\overline X,\overline Y) = \overline X^n - \overline Y^n - (\overline X -
\overline Y)^n$, where $n\geq 0$. 

It follows that $f(\overline X,\overline Y) \in \Ker(\overline d) \cap 
(\overline X\ \overline Y)$ iff $f$ is a linear span of the 
$ \overline X \big( \overline X^n - \overline Y^n - (\overline X -
\overline Y)^n \big) / (\overline X - \overline Y)$, $n\geq 0$ and 
belongs to $(\overline X\ \overline Y)$. This means that $f$ is a linear span 
of the same elements, where $n\geq 2$.  
\hfill \qed \medskip 

Let us now prove Theorem \ref{main:thm}. Let $\psi\in\t_3$ 
be a solution of (\ref{hex}) and (\ref{pent}), 
homogeneous of degree $n$. One checks that if $n=1$, then $\psi=0$; 
let us assume that $n\geq 2$. Recall that $\f_2\subset \t_3$ is 
the Lie subalgebra generated by $X = t_{13}$ and $Y = t_{23}$, and that  
$\t_3 = \kk \cdot (t_{12}+t_{13}+t_{23}) \oplus \f_2$. Since this is a graded
decomposition, we have $\psi\in\f_2'$, where $\f_2' = [\f_2,\f_2]$
is the degree $\geq 2$ part of $\f_2$ (it coincides with $\p$ defined in the
Introduction, since it coincides with $\t_3' = [\t_3,\t_3]$). 

Let us set $P_{k\ell} = \on{ad}(X)^{k-1}\on{ad}(Y)^{\ell-1}([X,Y])$
(here $k,\ell\geq 1$). 
Then $\p /[\p,\p]$ is an abelian Lie algebra with basis $[P_{k\ell}]$, 
$k,\ell\geq 1$. We have therefore $\psi = \sum_{k,\ell\geq 1, k+\ell = n}
a_{k\ell} P_{k\ell} + \psi'$, where $\psi'\in [\p,\p]$. We set 
$$
a(\overline X,\overline Y) := \sum_{k,\ell\geq 1,k+\ell = n}
a_{k\ell} \overline X^k\overline Y^\ell \in \kk[\overline X,\overline Y]. 
$$ 

\begin{lemma} The image of $P_{k\ell}$ in $F_2/I_2 \simeq \kk[\overline
X,\overline Y]$ is $(-1)^k \overline X^k \overline Y^\ell$. 
We have $[\p,\p] \subset I_2$. 
\end{lemma}

{\em Proof.} The first statmeent follows from the expansion of $P_{k\ell}$. 
Let us denote by $(F_2)_{>0} \subset F_2$ the subspace of all 
elements of positive valuation both in $X$ and in $Y$. 
Then $\p\subset (F_2)_{>0}$. So $[\p,\p] \subset 
((F_2)_{>0})^2 \subset I_2$, which proves the second statement. 
\hfill \qed \medskip 

Let us denote by $\overline\psi$ the image of $\psi$ in $F_2/I_2 \simeq
\kk[\overline X,\overline Y]$. 
Then $\overline\psi = \sum_{k,\ell\geq 1,k+\ell=n} (-1)^k a_{k\ell}
\overline X^k \overline Y^\ell = a(-\overline X,\overline Y)$. 

The image of $(F_2)_{>0}$ by the projection map $F_2 \to F_2/I_2 = 
\kk[\overline X,\overline Y]$ is the ideal 
$(\overline X\ \overline Y)$. Since $\psi\in \p$, we have 
$\overline\psi\in  (\overline X\ \overline Y)$. 

On the other hand, we have $\overline d(\overline\psi) = 0$. 
It then follows from Lemma \ref{ker:bar:d} that for some 
$\lambda\in\kk$, we have 
$a(-\overline X,\overline Y) = \lambda \overline X 
\big( \overline X^n - \overline Y^n - (\overline X - \overline Y)^n\big) 
/ (\overline X - \overline Y)$, i.e., 
$$
a(\overline X,\overline Y) = (-1)^{n+1}\lambda \overline X
{{ (\overline X + \overline Y)^n
- \overline X^n + (-1)^n \overline Y^n}\over{\overline X+\overline Y}}. 
$$

Recall that $p_{k\ell} = \on{ad}(A)^{k-1} \on{ad}(B)^{\ell-1}([A,B])$, 
where $A = t_{12}$ and $B = t_{23}$. 
Let $b_{k\ell}$ ($k,\ell>0$, $k+\ell=n$) be the coefficients such that 
$\psi\in \sum_{k,\ell\geq 1,k+\ell=n} b_{k\ell} p_{k\ell} + [\p,\p]$. We
set $b(\overline A,\overline B) = \sum_{k,\ell\geq 1,k+\ell = n} b_{k\ell} 
\overline A^k \overline B^\ell \in \kk[\overline A,\overline B]$.
Then $b(\overline A,\overline B)$ is the image of the class $[\psi]$
of $\psi$ in $\p/[\p,\p]$ under $i : \p/[\p,\p] \simeq 
(\overline A\ \overline B)$ defined in the Introduction. 

In general, the polynomials $a(\overline X,\overline Y)$ and 
$b(\overline A,\overline B)$ are related by  
$$
b(\overline A,\overline B) = 
-{{\overline A}\over {\overline A+\overline B}} 
a(-\overline A-\overline B,\overline B), 
$$ 
so in our case 
$$
b(\overline A,\overline B) = \lambda \big( \overline A^n + \overline B^n 
- (\overline A+\overline B)^n \big).  
$$

Now the image of condition (\ref{hex}) in $\p/[\p,\p]$ is that 
$b(\overline A,\overline B)$ satisfies 
$$
{{b(\overline A,\overline B)}\over{\overline A\ \overline B}}
+ {{b(\overline B,\overline C)}\over{\overline B\ \overline C}}
+ {{b(\overline C,\overline A)}\over{\overline C\ \overline A}} = 0,  
$$
where $\overline C = - \overline A - \overline B$. 

Now $\overline C \psi(\overline A,\overline B)
+ \overline A \psi(\overline B,\overline C)
+ \overline B \psi(\overline C,\overline A)
= - \lambda (1+(-1)^n)(\overline A^{n+1} 
+ \overline B^{n+1} + \overline C^{n+1})$. 

It follows that if $n$ is even, then the image of (\ref{hex}) 
implies $\lambda = 0$, 
therefore the image $[\psi]$ of $\psi$ in $\p/[\p,\p]$ is zero, 
and that if $n$ is odd, then the image of (\ref{hex}) is automatically
satisfied, so that $b(\overline A,\overline B)$ is proportional to 
$(\overline A + \overline B)^n - \overline A^n - \overline B^n$, 
i.e., $[\psi]$ is proportional to $[\sigma_n]$. This ends the proof 
of Theorem \ref{main:thm}. 
\hfill \qed \medskip 


\section{Proof of Corollary \ref{cor:Gamma}}

Recall that $\Assoc(\kk)$ is a torsor under the right action of 
a group $\on{GRT}(\kk)$. We will first prove: 

\begin{proposition} \label{prop:subtorsor}
Set $\Assoc^*(\kk) = \{\Phi\in \Assoc(\kk) | (\ref{Phi:Gamma})$ holds$\}$. 
Then $\Assoc^*(\kk)$ is stable under the action of $\on{GRT}(\kk)$
on $\Assoc(\kk)$. Therefore  
$\Assoc^*(\kk)$ is either $\emptyset$ or $\Assoc(\kk)$. 
\end{proposition}

{\em Proof of Proposition \ref{prop:subtorsor}.}
Let $\Phi\in\Assoc^*(\kk)$ and let $g\in\on{GRT}(\kk)$. 
We should prove that $\Phi * g$ satisfies (\ref{Phi:Gamma}). 

Recall that $\on{GRT}(\kk)$ is the semidirect product 
$\on{GRT}_1(\kk) \rtimes \kk^\times$, where $\on{GRT}_1(\kk)$
is the prounipotent group exponentiating $\grt_1(\kk)$, and
the action of $\kk^\times$ on $\on{GRT}_1(\kk)$ is the exponential of 
its action on $\grt_1(\kk)$ induced by the grading. 
So it suffices to check that $\Phi * g$ satisfies (\ref{Phi:Gamma})
when $g \in \kk^\times$, and when $g \in \on{GRT}_1(\kk)$. 

If $g = \mu\in \kk^\times$, then $\Phi * \mu = \Phi(\mu A,\mu B)$, 
therefore $\Phi * \mu$ satisfies (\ref{Phi:Gamma}) with 
$\zeta_{\Phi * \mu}(n) = \mu^n \zeta_\Phi(n)$. 

If $g\in \on{GRT}_1(\kk)$, then $g = \on{exp}(\psi)$, where 
$\psi\in \grt_1(\kk)$. We set $\Phi_t := \Phi * \exp(t\psi)$. 
According to Theorem \ref{main:thm}, there exist scalars 
$\mu_n\in \kk$ ($n$ odd $\geq 3$) such that $[\psi] = \sum_{n\ \on{odd},
n\geq 3} \mu_n [\sigma_n]$. Since $\psi\in\f_2(A,B)$, this means that 
$(\psi_B B)^{\on{ab}} = \sum_{n\ \on{odd}, n\geq 3} \mu_n
((\overline A + \overline B)^n - \overline A^n - \overline B^n)$. 

Let $\eps$ be a formal variable with $\eps^2=0$. Then 
$\Phi_{t+\eps} = \Phi_t + \eps \big( \Phi_t \psi + D_\psi(\Phi_t)\big)$, 
where $D_\psi$ is the derivation of $\f_2(A,B)$ such that 
$D_\psi(A) = [\psi,A]$, $D_\psi(B)=0$. 

Using the decompositions $\psi = \psi_A A + \psi_B B$, $\Phi_t = 
1 + (\Phi_t)_A A + (\Phi_t)_B B$, we get 
\begin{align*}
& \Phi_{t+\eps} = 1 + (\Phi_t)_A A + (\Phi_t)_B B  
\\ & + \eps \Big( \Phi_t \psi_A A + \Phi_t \psi_B B + D_\psi((\Phi_t)_A)A
+ D_\psi((\Phi_t)_B)B 
+ (\Phi_t)_A \big( \psi A - A (\psi_A A + \psi_B B)\big) \Big),  
\end{align*}
so 
$$
(\Phi_{t+\eps})_B B  = (\Phi_t)_B B + 
\eps \big( \Phi_t \psi_B B + D_\psi((\Phi_t)_B)B 
- (\Phi_t)_A A  \psi_B B \big).   
$$
Let us apply the abelianization to this formula. 
Since $\Phi_t\in \on{exp}(\wh\f_2(A,B))$, we have 
$\Phi_t^{\on{ab}} = 1$ and so 
$((\Phi_t)_A A + (\Phi_t)_B B)^{\on{ab}} = 0$. 
Therefore 
$$
(d/dt)(((\Phi_t)_B B)^{\on{ab}}) = (\psi_B B)^{\on{ab}}
\big( 1- ((\Phi_t)_A A)^{\on{ab}} \big) 
= (\psi_B B)^{\on{ab}}
\big( 1 + ((\Phi_t)_B B)^{\on{ab}} \big) . 
$$
Therefore $1+((\Phi_t)_B B)^{\on{ab}} = 
\big( 1+(\Phi_B B)^{\on{ab}} \big) \exp\big( t 
(\psi_B B)^{\on{ab}}
\big)$, and with $t=1$ this gives 
$1+((\Phi * g)_B B)^{\on{ab}} = 
\big( 1+(\Phi_B B)^{\on{ab}} \big) \exp\big(  
(\psi_B B)^{\on{ab}} \big)$. 

Since $\Phi$ satisfies (\ref{Phi:Gamma}), we get 
$1+((\Phi * g)_B B)^{\on{ab}} = 
\Gamma_{\Phi * g}(\overline A + \overline B) / 
(\Gamma_{\Phi * g}(\overline A)  
\Gamma_{\Phi * g}(\overline B) )$, where 
$\Gamma_{\Phi * g}(s) = \Gamma_\Phi(s) 
\exp(\sum_{n\ \on{odd},n\geq 3} \mu_n s^n)$, i.e. 
$\Phi * g$ satisfies (\ref{Phi:Gamma}) with 
$\zeta_{\Phi * g}(n) = \zeta_\Phi(n)  - n \mu_n$
for $n$ odd $\geq 3$, and 
$\zeta_{\Phi * g}(n) = \zeta_\Phi(n)$ for $n$ even $\geq 2$. 
\hfill \qed \medskip 

Let us now prove Corollary \ref{cor:Gamma}. 
Proposition \ref{prop:subtorsor} implies that 
$\Assoc^*(\kk)$ is either $\emptyset$ or $\Assoc(\kk)$. 

Let $\kk$ and $\kk'$ be fields of characteristic $0$. 
It is immediate that if $\kk \subset \kk'$ and 
$\Assoc^*(\kk') = \Assoc(\kk')$, then 
$\Assoc^*(\kk) = \Assoc(\kk)$. 
On the other hand, if $\kk \subset \kk'$ and  
$\Assoc^*(\kk) = \Assoc(\kk)$, then 
$\Assoc^*(\kk') = \Assoc(\kk')$: indeed, 
Proposition 5.3 of \cite{Dr} implies that 
$\Assoc(\kk)\neq \emptyset$, so 
$\Assoc^*(\kk) \neq \emptyset$; 
we have obviously $\Assoc^*(\kk) \subset \Assoc^*(\kk')$, 
hence $\Assoc^*(\kk') \neq \emptyset$; then Proposition 
\ref{prop:subtorsor} implies that $\Assoc^*(\kk') = \Assoc(\kk')$. 
It follows that if for some $\kk$, $\Assoc^*(\kk) \neq \emptyset$, 
then $\Assoc^*(\kk) = \Assoc(\kk)$ for any $\kk$. We will now prove that 
$\Assoc^*(\CC) \neq \emptyset$. 

Let $\Phi_{\on{KZ}}$ be the Knizhnik-Zamolodchikov associator defined as in 
\cite{Dr} as the renormalized holonomy from $0$ to $1$ of the differential
equation $G'(z) = ({A\over z} + {B\over{z-1}})G(z)$. 
Then $(2\pi i,\Phi_{\on{KZ}}) \in \Assoc(\CC)$
satisfies (\ref{Phi:Gamma}) with $\zeta_\Phi(n) = \zeta(n)$ 
for any $n\geq 2$. Indeed, in \cite{Dr}, (2.15), it is proved 
that 
$$
[\log\Phi_{\on{KZ}}] = \exp \big( \sum_{n\geq 2} {{\zeta(n)}\over n}
(\overline A^n + \overline B^n - (\overline A + \overline B)^n)\big) - 1. 
$$
Then $(\Phi_{\on{KZ}})_A = {{\Phi_{\on{KZ}} - 1}\over
{\log\Phi_{\on{KZ}}}} (\log\Phi_{\on{KZ}})_A$, 
$(\Phi_{\on{KZ}})_B = {{\Phi_{\on{KZ}} - 1}\over
{\log\Phi_{\on{KZ}}}} (\log\Phi_{\on{KZ}})_B$, therefore 
$(\Phi_{\on{KZ}})_B^{\on{ab}} = (\log\Phi_{\on{KZ}})_B^{\on{ab}}
= [\log\Phi_{\on{KZ}}] / \overline B$ (the last equality follows 
from $\log\Phi_{\on{KZ}}\in\p$). So 
$$
1 + ((\Phi_{\on{KZ}})_B B)^{\on{ab}} = 
{{\Gamma_{\on{mod}}(\overline A + \overline B)}
\over{\Gamma_{\on{mod}}(\overline A)\Gamma_{\on{mod}}(\overline B)}}, 
$$
where $\Gamma_{\on{mod}}(u) = \exp(\sum_{n\geq 2} - {{\zeta(n)}\over n} u^n)$
is related to the $\Gamma$-function by 
$\Gamma_{\on{mod}}(u) = e^{\gamma u} / (-u\Gamma(-u))$, where 
$\gamma$ is the Euler-Mascheroni constant. It follows that 
$(\Phi_{\on{KZ}},2\pi i) \in \Assoc^*(\CC)$, therefore for any 
$\kk$, $\Assoc^*(\kk) = \Assoc(\kk)$. \hfill \qed 

\end{document}